\documentstyle[amssymb,amsfonts]{amsart}

\def\R{{\hbox{\bf R}}}

\def\P{{\hbox{\bf P}}}
\def\E{{\hbox{\bf E}}}

\font \roman = cmr10 at 10 true pt

\def\be#1{ \begin{equation}\label{#1} }

\def\bas{\begin{align*}}
\def\eas{\end{align*}}
\def\bi{\begin{itemize}}
\def\ei{\end{itemize}}

\def\dist{{\hbox{\roman dist}}}

\def\eps{\varepsilon}

\def \endprf{\hfill  {\vrule height6pt width6pt depth0pt}\medskip}
\def\emph#1{{\it #1}}
\def\textbf#1{{\bf #1}}

\def\Span{\hbox{ \rm Span} \,\,\, }

\theoremstyle{plain}
  \newtheorem{theorem}[subsection]{Theorem}
  \newtheorem{conjecture}[subsection]{Conjecture}

  \newtheorem{lemma}[subsection]{Lemma}

\theoremstyle{remark}
  \newtheorem{remark}[subsection]{Remark}

\theoremstyle{definition}

\include{psfig}

\begin{document}

\title[Littlewood-Offord in high dimensions]{The Littlewood-Offord problem in high dimensions
and a conjecture of Frankl and F\"uredi}

\author{Terence Tao}
\address{Department of Mathematics, UCLA, Los Angeles CA 90095-1555}
\email{tao@@math.ucla.edu}
\thanks{T. Tao is supported by NSF Research Award DMS-0649473, the NSF Waterman award and a grant from the MacArthur Foundation.}

\author{Van Vu}
\address{Department of Mathematics, Rutgers, New Jersey, NJ 08854}
\email{vanvu@@math.rutgers.edu}
\thanks{V. Vu is supported by research grants DMS-0901216 and AFOSAR-FA-9550-09-1-0167.}

\begin{abstract}
We give a new  bound on the probability that the random sum 
$\xi_1 v_1 +\dots + \xi_n v_n$ belongs to a ball of fixed radius, where the $\xi_i$ are iid Bernoulli random variables and the $v_i$ are vectors in $\R^d$. As an application, we 
prove a conjecture of  Frankl and F\"uredi (raised in 1988), which can be seen as the high dimensional version of the classical Littlewood-Offord-Erd\H os theorem. 
\end{abstract}

\maketitle

\section{Introduction}

Let $V=\{v_1, \dots, v_n\}$ be a (multi-)set of $n$ vectors in $\R^d$. Consider the random sum

$$ X_{V} := \xi_1 v_1 + \dots \xi_n v_n $$

\noindent where $\xi_i$ are i.i.d. Bernoulli random variables (each $\xi_i$ takes values $1$ and $-1$ with probability $1/2$ each). 

\noindent The famous \emph{Littlewood-Offord problem} (posed in 1943 \cite{LO})  is to  estimate the \emph{small ball probability}

$$p_d(n, \Delta)= \sup_{V,B} \P( X_V \in B )$$

\noindent where the supremum is taken over all multi-sets $V = \{v_1,\ldots,v_n\}$ of $n$ vectors of length at least one and all
closed balls $B$ of radius $\Delta$ (this problem is also sometimes referred to as the \emph{small ball problem} in the literature).  
Here and later, $d$ and $\Delta$ are fixed. The asymptotic notation $X=O(Y)$ or  (equivalently) $X \ll Y$ will be used with the assumption that $n$ tends to infinity; thus the implied constant in the $O()$ notation can depend on $d$ and $\Delta$ but not on $n$. 

The more combinatorial (but absolutely equivalent) way to pose the problem is to ask for the maximum
number of subsums of $V$ falling into a  ball of radius $\Delta/2$. We prefer the probabilistic setting as it  more convenient and easier to generalize.

Shortly after the paper of Littlewood-Offord, Erd\H os \cite{Erd1} determined $p_1 (n, \Delta)$, solving the problem completely in one dimension.  Define $s:=\lfloor \Delta \rfloor +1$.

\begin{theorem}[Erd\H os' Littlewood-Offord inequality] \label{Erd1} Let $S(n,m)$ denote the sum of the largest $m$ binomial coefficients
${n \choose i}, 0 \le i \le n$. Then

$$p_1(n,\Delta) = 2^{-n} S(n, s). $$
\end{theorem}


The situation for higher dimension is  more complicated, and 
there has been a series of papers devoted to its study 
(see  \cite{Kat, Kle1, Kle2, Kle3, GLOS, Hal, FF, S1, S2} and the references therein). In particular,    Frankl and F\"uredi \cite{FF}, sharpening several earlier results, proved

\begin{theorem}[Frankl-F\"uredi's Littlewood-Offord inequality]\label{FF1} For any fixed $d$ and $\Delta$
\begin{equation}\label{panda}
p_d(n,\Delta) = (1+o(1)) 2^{-n} S(n,s). 
\end{equation} \end{theorem}

This result is asymptotic. In view of Theorem \ref{Erd1}, it is natural to ask if 
one can have the exact estimate \begin{equation} \label{panda-2} p_d (n, \Delta) = 2^{-n} S(n,s), \end{equation}  
 which can be seen as the 
high dimensional  generalization of Erd\H os' result.  However, it  has turned out that in general this is not true.  It was observed in \cite{Kle2, FF}  that \eqref{panda-2}  fails if $s \geq 2$ and
\begin{equation}\label{kappas}
 \Delta   > \sqrt{(s-1)^2+1} .
\end{equation}

Take $v_1= \dots =v_{n-1} =e_1$ and $v_n=e_2$, where $e_1, e_2$ are two orthogonal unit vectors. For this system, there is a ball $B$ of radius $\Delta$ such that $\P(X_V \in B) > S(n,s)$. 

Frankl and F\"uredi conjectured (\cite[Conjecture 5.2]{FF})

\begin{conjecture} \label{conj:FF}
Let $\Delta,d$ be fixed. If  $s-1 \le  \Delta < \sqrt {(s-1)^2 +1}$ and 
$n$ is sufficiently large, then 

$$p_d(n, \Delta) = 2^{-n} S(n,s). $$
\end{conjecture} 

The conjecture has been confirmed for $s=1$ by an important result of Kleitman \cite{Kle1} and for $s=2,3$ 
by Frankl and F\"uredi \cite{FF}
(see the discussion prior to \cite[Conjecture 5.2]{FF}). For all other cases, the conjecture has been open. On the other hand, Frankl and F\"uredi  showed that 
\eqref{panda-2} holds under a stronger assumption that 
$s-1 \le  \Delta \le  (s-1) + \frac{1}{10s^2}$.

\vskip2mm

In this short paper,  we first prove the following general estimate:

\begin{theorem} \label{TV1} Let $ V=\{ v_1, \dots, v_n\}$ be a multi-set of  vectors in $\R^d$ with the property that for any hyperplane $H$, one has $\dist(v_i,H) \geq 1$ for at least $k$ values of $i=1,\ldots,n$. Then for any unit ball $B$, one has
$$\P (X_V \in B) = O(k^{-d/2}). $$ The hidden constant in the $O()$ notation here depends on $d$, but not on $k$ and $n$. 
\end{theorem}


As an application, we prove Conjecture \ref{conj:FF} in full generality and also give a new proof 
for Theorem \ref{FF1}.  This will be done in the next section. The remaining two sections are devoted to the proof of Theorem \ref{TV1}. 
 
\section{Proof of Theorem \ref{FF1} and Conjecture \ref{conj:FF} }\label{ftsec}

We now assume Theorem \ref{TV1} is true, and use it to first prove Theorem \ref{FF1}. 
We will  induct on the dimension $d$.  The case $d=1$ follows from Theorem \ref{Erd1}, so we assume that $d \geq 2$ and that the claim has already been proven for smaller values of $d$.  The lower bound
$$p_d(n,\Delta) \geq p_1(n,\Delta) =  2^{-n} S(n,s)$$
is clear, so it suffices to prove the upper bound
$$p_d(n,\Delta) \leq (1+o(1)) 2^{-n} S(n,s).$$

Fix $\Delta$, and let $\eps > 0$ be a small parameter to be chosen later.  Suppose the claim failed, then there exists $\Delta > 0$ such that for arbitrarily large $n$, there exist a family $V = \{v_1,\ldots,v_n\}$ of vectors in $\R^d$ of length at least $1$ and a  ball $B$ of radius $\Delta$ such that
\begin{equation}\label{xv}
\P( X_V \in B ) \geq (1+\eps) 2^{-n} S(n,s).
\end{equation}
In particular, from Stirling's approximation one has
$$ \P(X_V \in B) \gg n^{-1/2}.$$

Assume $n$ is sufficiently large depending on $d,\eps$, and that $V, B$ is of the above form.  Applying the pigeonhole principle, we can find a ball $B'$ of radius $\frac{1}{\log n}$ such that
$$ \P(X_V \in B') \gg n^{-1/2} \log^{-d} n.$$

Set $k := n^{2/3}$.  Since $d \ge 2$ and $n$ is large, we have
$$ \P(X_V \in B') \geq C k^{-d/2}$$
for any fixed constant $C$.  Applying Theorem \ref{TV1} in the contrapositive (rescaling by $\log n$), we conclude that there exists a hyperplane $H$ such that $\dist(v_i,H) \leq 1/\log n$ for at least $n-k$ values of $i=1,\ldots,n$.

Let $V'$ denote the orthogonal projection to $H$ of the vectors $v_i$ with $\dist(v_i,H) \leq 1/\log n$.  By conditioning on the signs of all the $\xi_i$ with $\dist(v_i,H) > 1/\log n$, and then projecting the sum $X_V$ onto $H$, we conclude from \eqref{xv} the existence of a $d-1$-dimensional ball $B'$ in $H$ of radius $\Delta$ such that
$$ \P( X_{V'} \in B' ) \geq (1+\eps) 2^{-n} S(n,s).$$
On the other hand, the vectors in $V'$ have magnitude at least $1-1/\log n$.  If  $n$ is sufficiently large depending on $d,\eps$ this contradicts the induction hypothesis (after rescaling the $V'$ by $1/(1-1/\log n)$ and identifying $H$ with $\R^{n-1}$ in some fashion).  This concludes the proof of \eqref{panda}.

Now we turn to the proof of Conjecture \ref{conj:FF}. 
 We can assume $s \ge 3$, as the remaining cases have already been treated. 
   If the conjecture failed, then there exist arbitrarily large $n$ for which there exist a family $V = \{v_1,\ldots,v_n\}$ of vectors in $\R^d$ of length at least $1$ and a  ball $B$ of radius $\Delta$  such that
\begin{equation}\label{xv-2}
\P( X_V \in B ) > 2^{-n} S(n,s).
\end{equation}

 By iterating the argument used to prove \eqref{panda}, we may find a one-dimensional subspace $L$ of $\R^d$ such that $\dist(v_i,L) \ll 1/\log n$ for at least $n-O(n^{2/3})$ values of $i=1,\ldots,n$.  By reordering, we may assume that $\dist(v_i,L) \ll 1/\log n$ for all $1 \leq i \leq n-k$, where $k = O(n^{2/3})$.

Let $\pi: \R^d \to L$ be the orthogonal projection onto $L$.  We divide into two cases.  
The first case is when  $|\pi(v_i)| > \frac{\Delta}{s}$ for all $1 \leq i \leq n$.  We then use the trivial bound
$$\P( X_V \in B ) \leq \P( X_{\pi(V)} \in \pi(B) ).$$
If we rescale Theorem \ref{Erd1} by a factor slightly less than $s/\Delta$, we see that
$$ \P( X_{\pi(V)} \in \pi(B )) \leq 2^{-n} S(n,s) $$
which contradicts \eqref{xv-2}.

In the second case, we assume $|\pi(v_n)| \leq \Delta/s$.  We let $V'$ be the vectors $v_1,\ldots,v_{n-k}$, then by conditioning on the $\xi_{n-k+1},\ldots,\xi_{n-1}$ we conclude the existence of a unit ball $B'$ such that
$$ \P(X_{V'} + \xi_n v_n  \in B') \geq \P( X_V \in B ).$$

Let $x_{B'}$ be the center of $B'$.  Observe  that 
if $X_{V'} + \xi_n v_n \in B'$ (for any value of $\xi_n$) then  $|X_{\pi(V')} - \pi(x_{B'})| \leq \Delta + \frac{\Delta}{s}$.  Furthermore, if  $|X_{\pi(V')} - \pi(x_{B'})| > \sqrt{\Delta^2-1}$, then the parallelogram law shows that $X_{V'} + v_n$ and $X_{V'}-v_n$ cannot both lie in $B'$, and so conditioned on  $|X_{\pi(V')} - \pi(x_{B'})| > \sqrt{\Delta^2-1}$, the probability that $X_{V'} + \xi_n v_n \in B'$ is at most $1/2$.

 We conclude that
\begin{align*}
& \P(X_{V'} + \xi_n v_n  \in B') \\ &\leq  \P( |X_{\pi(V')} - \pi(x_{B'})| \leq \sqrt{\Delta^2-1} )
+ \frac{1}{2}  \P( \sqrt {\Delta^2-1} < |X_{\pi(V')} - \pi(x_{B'})| \leq \Delta+\frac{\Delta}{s} ) \\
&=\frac{1}{2} \Big(\P( |X_{\pi(V')} - \pi(x_{B'})| \leq \sqrt{\Delta^2-1} )+ \P(|X_{\pi(V')} - \pi(x_{B'})| \leq \Delta+\frac{\Delta}{s} ) \Big).
\end{align*}

However, note that all the elements of $\pi(V')$ have magnitude at least $1-1/\log n$.  
Assume, for a moment,  that $\Delta$ satisfies 
\begin{equation}\label{deltas}
\sqrt{\Delta^2-1} < s-1 \le \Delta <   \Delta +\frac{\Delta}{s}  < s.
\end{equation}
  From Theorem \ref{Erd1} (rescaled by $(1-1/\log n)^{-1}$), we conclude  that
$$
\P( |X_{\pi(V')} - \pi(x_{B'})| \leq \sqrt{\Delta^2-1} ) \leq 2^{-(n-k)} S(n-k,s-1)$$
and
$$
\P( |\pi(X_{V'}) - \pi(x_{B'})| \leq \Delta+\frac{\Delta}{s} ) \leq 2^{-(n-k)} S(n-k,s).$$

On the other hand, by Stirling's formula (if $n$ is sufficiently large) we have
$$ \frac{1}{2}  (2^{-(n-k)} S(n-k,s-1)) + \frac{1}{2} 2^{-(n-k)} S(n-k,s) = \sqrt {\frac{2}{\pi}} \frac{s-1/2+o(1)}{n^{1/2}}$$
while
$$ 2^{-n} S(n,s) = \sqrt {\frac{2}{\pi}} \frac{s+o(1)}{n^{1/2}}$$
and so we contradict \eqref{xv-2}. 

\vskip2mm

  An inspection of the above argument shows that all we need on $\Delta$ are the conditions \eqref{deltas}. To satisfy the first inequality in \eqref{deltas}, we need
  $\Delta < \sqrt {(s-1)^2 +1}$.  Moreover, once $s-1 \le \Delta < \sqrt{(s-1)^2+1}$, 
  one can easily check that  $\Delta +\frac{\Delta}{s} < s$
holds automatically for any $s \ge 3$,  concluding  the  proof.

\section{Proof of Theorem \ref{TV1}}

Let $d,n,k,V$ be as in Theorem \ref{TV1}.  We allow all implied constants to depend on $d$.

By Ess\'een's concentration inequality (see \cite{Hal}, \cite{TVinv}, or \cite[Lemma 7.17]{TVbook}), we have for any unit ball $B$ that
$$\P (X_V \in B) \ll  \int_{\zeta \in \R^d: |\zeta| \leq 1} |\E( e( \zeta \cdot X_V ) )|\ d\zeta.$$
and $e(x) := e^{2\pi \sqrt{-1} x}$.  
From the definition of $X_V$ and independence we have
$$ \E( e( \zeta \cdot X_V ) ) = \prod_{j=1}^n \E (e(\zeta \cdot \xi_j v_j)) =  \prod_{j=1}^n \cos( \pi \zeta \cdot v_j ).$$
Denoting by $\|\theta\|$  the distance from $\theta$ to the nearest integer and using the elementary bound 
$ |\cos(\pi \theta)| \leq \exp( - \frac{ \| \theta \|^2 }{100} )$ (whose proof is left as an exercise),
we reduce to showing the bound
\begin{equation} \label{eqn:main} Q\ll k^{-d/2}. \end{equation} 
where
\begin{equation}\label{zd}
Q := \int_{\zeta \in \R^d: |\zeta| \leq 1} \exp( - \frac{1}{100} \sum_{v \in V} \| \zeta\cdot v\|^2 )\ d\zeta.
\end{equation}

To show \eqref{zd}, our main technical tool  is the following lemma, whose proof is deferred to the next  section. 

\begin{lemma}\label{dspan} Let $w_1,\dots,w_d \in \R^d$ be such that
$\dist (w_j, \Span \{ w_1,\dots,  w_{j-1} \} ) \ge 1$ for  each $1 \leq j \leq d$, where $\Span\{w_1,\ldots,w_{j-1}\}$ is the linear span of the $w_1,\ldots,w_{j-1}$, and $\dist$ denotes Euclidean distance. Then for any $\lambda > 0$,
$$ \int_{\zeta \in \R^d: |\zeta| \leq 1} \exp( - \lambda \sum_{j=1}^d \|\zeta \cdot w_j \|^2 )\ d\zeta = O( (1 + \lambda)^{-d/2} ).$$
\end{lemma}

With this lemma  in hand, we conclude the proof as follows.  By shrinking $k$, we may assume that $k=dl$ for some integer $l$. Let $v_{0,1},\ldots,v_{0,l}$ be $l$ elements of $V$, and let $V_1 := V \backslash \{v_{0,1},\ldots,v_{0,l}\}$.  Then we can write
$$
Q = \int_{\zeta \in \R^d: |\zeta| \leq 1} \exp( - \frac{1}{100} \sum_{v \in V_1} \| \zeta\cdot v\|^2 )
\prod_{j=1}^l \exp( - \frac{1}{100} \|\zeta \cdot v_{0,j} \|^2 )\ d\zeta.$$
Applying H\"older's inequality, we conclude the existence of a $j=1,\ldots,l$ such that
$$
Q \leq \int_{\zeta \in \R^d: |\zeta| \leq 1} \exp( - \frac{1}{100} \sum_{v \in V_1} \| \zeta\cdot v\|^2 )
\exp( - \frac{l}{100} \|\zeta \cdot v_{0,j} \|^2 )\ d\zeta.$$
Write $w_1 := v_{0,j}$.  If $d=1$, we stop at this point.  Otherwise, we choose $l$ elements $v_{1,1},\ldots,v_{1,l}$ be $l$ elements of $V_1$ which lie at a distance at least $1$ from the span $\Span\{w_1\}$ of $w_1$; such elements can be found thanks to the hypotheses of Theorem \ref{TV1}.  We write $V_2 := V_1 \backslash \{v_{1,1},\ldots,v_{1,l}\}$.  By using H\"older's inequality as before, we can find $j=1,\ldots,l$ such that
$$
Q \leq \int_{\zeta \in \R^d: |\zeta| \leq 1} \exp( - \frac{1}{100} \sum_{v \in V_2} \| \zeta\cdot v\|^2 )
\exp( - \frac{l}{100} \|\zeta \cdot w_1 \|^2 ) \exp( - \frac{l}{100} \|\zeta \cdot v_{1,j} \|^2 )
\ d\zeta.$$
We then set $w_2 := v_{1,j}$.  We repeat this procedure $d-1$ times, eventually obtaining
$$
Q \leq \int_{\zeta \in \R^d: |\zeta| \leq 1} \exp( - \frac{1}{100} \sum_{v \in V_d} \| \zeta\cdot v\|^2 )
\exp( - \frac{l}{100} \sum_{i=1}^d \|\zeta \cdot w_i \|^2 )
\ d\zeta$$
for some $w_1,\ldots,w_d$ with the property that $\dist(w_i, \Span\{w_1,\ldots,w_{i-1}\}) \geq 1$ for all $1 \leq i \leq d$, and where $V_d$ is a subset of $V$ of cardinality at least $n-k$.  If we then trivially bound
$\exp( - \frac{1}{100} \sum_{v \in V_d} \| \zeta\cdot v\|^2 )$ by one, the claim follows from Lemma \ref{dspan}.

\begin{remark} An inspection of the argument reveals that Theorem \ref{TV1} still holds if one replaces the Bernoulli random variables by 
more general ones. For example, it suffices to assume that  $\xi_1, \dots, x_n$ are independent random variables satisfying 
$| \E e(x_i t )|  \le (1-\mu) + \mu \cos \pi t$ for any real number $t$, where $0 < \mu \le 1$ is a constant. Indeed, with this assumption we have
$$| \E e(x_i t )|  \leq \exp(-c_\mu \|t\|^2)$$
for all $t$ and some $c_\mu > 0$, and the rest of the argument can then be continued with $c_\mu$ playing the role of the constant $1/100$.

It is easy to see that
if there are constants $K, \epsilon$ such that  the support of 
every $\xi_i$ belongs to $\{-K, \dots, K\}$, and $\P(\xi =j)  \le 1-\epsilon$ for all $-K \le j \le K$, then all $\xi_i$ are $\mu$-bounded for some 
$0 < \mu \le 1$ depending on $K$ and $\epsilon$. 
\end{remark}

\section{Proof of Lemma \ref{dspan} }

The only remaining task is to show Lemma \ref{dspan}.  We are going to prove this lemma in the following, slightly more general but more convenient, form. 

\begin{lemma}\label{dspan2} Let $w_1,\dots,w_d \in \R^d$ be such that
$\dist (v_j, \Span \{ w_1,\dots w_{j-1} \} ) \ge 1$, for  each $1 \leq j \leq d$.  Let $u_1, \dots, u_d$ be arbitrary numbers. Then for any $\lambda > 0$,

\begin{equation} \label{bound1}  \int_{\zeta \in \R^d: |\zeta| \leq 1} \exp( - \lambda \sum_{j=1}^d \|\zeta \cdot w_j + u_j \|^2 )\ d\zeta \ll (1 + \lambda)^{-d/2}. \end{equation} 
Again, we allow all implied constants to depend on $d$.
\end{lemma}

We first consider the  case $d=1$. It this case the claim is equivalent to 

$$\int_{\zeta \in \R; |\zeta-u_1| \le w_1 } \exp( -\lambda \| \zeta\| ^2) d \zeta = O (\frac{|w_1|}{\sqrt {1+ \lambda}} ),$$

which follows from periodicity of the function $\| \zeta \|$  and the elementary estimate 

$$\int_{-1}^1 \exp( \| -\lambda \zeta \| ^2 ) d \zeta = O(\frac{1}{\sqrt {1+\lambda} } ), $$ whose proof is left as an exercise. 

To handle the general case, we use Fubini's theorem and induction on $d$. 
 By Gram-Schmidt orthogonalization, we can find an  orthonormal basis $\{e_1, \dots, e_d \}$ of $\R^d$. 
 such that  $\Span \{w_1, \dots, w_j \} = \Span \{e_1, \dots, e_j\}$, for all $1\le j \le d$. 
 Suppose that the desired claim holds for $d-1$. For a vector $\zeta \in \R^d$, write 
 
 $$\zeta := \zeta' + \zeta_d e_d $$ where $\zeta' \in \Span \{e_1, \dots, e_{d-1} \} $ and $\zeta_d \in \R$. 
 The left hand side of \eqref{bound1} can be rewritten as

 $$ \int_{|\zeta'| \le 1}  \Big[ \exp( - \lambda \sum_{j=1}^d \|\zeta \cdot w_j + u_j \|^2 ) \int_{|\zeta_d| \le 1} \exp\Big(-\lambda \| \zeta_d (e_d \cdot w_d) + (\zeta' \cdot w_d + u_d) \|^2 \Big) d \zeta_d \Big] d \zeta' . $$
 
 By the case $d=1$, the inner integral is $O(\frac{1}{\sqrt {\lambda+1}})$, uniformly in $\zeta'$. The claim now follows from the induction hypothesis.

\end{document}